# A note on percolation in cocycle measures

Ronald Meester[1]

*Vrije Universiteit, Amsterdam*

**Abstract:** We describe infinite clusters which arise in nearest-neighbour percolation for so-called cocycle measures on the square lattice. These measures arise naturally in the study of random transformations. We show that infinite clusters have a very specific form and direction. In concrete situations, this leads to a quick decision whether or not a certain cocycle measure percolates. We illustrate this with two examples which are interesting in their own right.

## 1. Introduction

Much of Mike's work in probability and percolation theory has been inspired by his background in ergodic theory. His ergodic-theoretical viewpoint of spatial stochastic models turned out to be very fruitful, both for answering long standing open questions as for generating new problems. Among many other things, Mike taught me how to think 'ergodically', and I have enjoyed the interplay between probability and ergodic theory ever since.

In this note, we will further illustrate this interplay in a concrete situation; we discuss some percolation properties of a particular class of random colourings of the nearest neighbour edges of the square lattice $\mathbb{Z}^2$. The (probability) measures in this note are related to measure-preserving random transformations, and for reasons that will become clear we shall call it the class of *cocycle measures*.

We consider colourings of the edges of $\mathbb{Z}^2$ with two colours, red and blue, with the *cocycle property*. This property can be reformulated as follows. Consider four nearest-neighbour edges forming a square. When you travel in two steps from southwest to north-east along this square, the number of blue and red edges you see along the way does not depend on the route you take. A second way of defining this class of colourings is as follows. Take any two vertices $x$ and $y \in \mathbb{Z}^2$, and consider a vertex-self-avoiding *path* $\pi$ between $x$ and $y$, i.e. a sequence of distinct edges $e_i = (u_i, v_i)$ $i = 1, \ldots, k$ such that $u_1 = x$, $v_k = y$, $v_i = u_{i+1}$ for $i = 1, \ldots, k-1$. When travelling along $\pi$ from $x$ to $y$, we travel edges horizontally to the right, vertically upwards, horizontally to the left or vertically downwards. We collect the first two types of edges in a set $\pi^+$, and the last two in a set $\pi^-$. Now consider the number of red edges in $\pi^+$, minus the number of red edges in $\pi^-$ and call this number $f_1(\pi)$. Similarly, $f_2(\pi)$ is defined as the number of blue edges in $\pi^+$ minus the number of blue edges in $\pi^-$. The requirement we impose on the configurations is that $(f_1(\pi), f_2(\pi))$ is the same for all paths $\pi$ from $x$ to $y$.

Motivation for this type of measures can for instance be found in Burton, Dajani and Meester (1998). Indeed, the last characterisation above is in fact a formulation of a so called cocycle-identity, but this will play no role in the present note. We give some examples of cocycle measures in the last section of this note.

---

[1]Deptartment of Mathematics, Vrije Universiteit Amsterdam, De Boelelaan 1081, 1081 HV Amsterdam, The Netherlands, e-mail: `rmeester@cs.vu.nl`

Let $\mu$ be a stationary, ergodic (with respect to the group of all translations of $\mathbb{Z}^2$) cocycle measure. We are interested in percolation properties of $\mu$, i.e. we are interested in the question whether or not infinite red or blue clusters exist, and if so, how many. It is easy to come up with examples in which both blue and red edges percolate; for instance take the measure $\mu$ which makes all horizontal edges blue and all vertical edges red. On the other hand, it is just as easy to find an example where neither the blue nor the red edges percolate; just colour the four edges of every second square blue and the remaining edges red, and choose the origin randomly as to get something stationary. We leave it to the reader to find an easy example of a measure for which exactly one colour percolates.

The goal of the present note is to discuss some general percolation properties for this type of measures, which in concrete examples lead to a quick decision whether or not a given measure actually percolates. For reasons that will become clear soon, we will no longer speak about red or blue edges, but about edges labelled 0 or 1, and from now on this refers to a number, not a colour. We shall concentrate on percolation of edges labelled 0. Of course, this is in some sense arbitrary, but 0's really seem to have advantages over 1's, as we shall see.

The next section gives some general background on cocycle measures. Section 3 deals with general facts about percolation in cocycle measures, and the last section is devoted to a number of examples. The first example in the last section was the motivation to study percolation properties of cocycle measures; in this example we answer a question which was asked by T. Hamachi.

## 2. General background

We start with some notation. The expection of the label of a horizontal edge is denoted by $h$, the expection of the label of a vertical edges by $v$. To avoid trivial situations, we assume that $0 < h, v < 1$. We write $f(z) = f(z_1, z_2) = f(z_1, z_2, \omega)$ for the sum of the labels in $\pi^+$ minus the sum of the labels in $\pi^-$, where $\pi$ is an arbitrary self-avoiding path from 0 to $z = (z_1, z_2)$. $L_1$ distance is denoted by $\|\cdot\|$.

Note that our weak assumption on ergodicity of $\mu$ does not imply that the right or up shift are individually ergodic. However, since by the defining property of cocycle measures we have that

$$|(f(n,m) - f(0,m)) - (f(n,m+1) - f(0,m+1))| \leq 2,$$

the limits $\lim_{n\to\infty}(f(n,m) - f(0,m))/n$, which exist by stationarity, are invariant under both horizontal and vertical translations and therefore a.s. constant. This constant has to be $h$ then. A similar remark is valid for vertical limits.

The cone
$$\{(x,y) \in \mathbb{Z}^2 : \alpha - \epsilon \leq \frac{y}{x} \leq \alpha + \epsilon\}$$
is denoted by $C(\alpha, \epsilon)$. Throughout, $\mu$ denotes a cocycle measure. The following lemma is taken from Dajani and Meester (2003).

**Lemma 2.1.** *Let $\{(k_n, m_n)\}$ be a sequence of vectors in $\mathbb{Z}^2$.*

(i) *Suppose that $(k_n, m_n) \to (c_1 \cdot \infty, c_2 \cdot \infty)$ for some $c_1, c_2 \in \{1, -1\}$ and in addition that $\frac{m_n}{k_n} \to \alpha \in [-\infty, \infty]$. Then*

$$\frac{f(k_n, m_n)}{|k_n| + |m_n|} \to \frac{c_1}{1 + |\alpha|} h + \frac{c_2 |\alpha|}{1 + |\alpha|} v$$

in $\mu$-probability as $n \to \infty$. (The quotient $\frac{1}{1+\infty}$ is to be interpreted as 0 and $\frac{\infty}{1+\infty}$ as 1.)

(ii) Suppose that $\{k_n\}$ is bounded and $m_n \to c_3 \cdot \infty$ for some $c_3 \in \{1, -1\}$. Then

$$\frac{f(k_n, m_n)}{|k_n| + |m_n|} \to c_3 v$$

in $\mu$-probability as $n \to \infty$.

(iii) Suppose that $\{m_n\}$ is bounded and $k_n \to c_4 \cdot \infty$ for some $c_4 \in \{1, -1\}$. Then

$$\frac{f(k_n, m_n)}{|k_n| + |m_n|} \to c_4 h$$

in $\mu$-probability as $n \to \infty$.

This leads to

**Lemma 2.2.** *Let $\alpha \in (-\infty, \infty)$. Then for any $\epsilon > 0$, there a.s. exist $N_\epsilon > 0$ and $\delta_\epsilon > 0$ such that whenever $|m_n|, |k_n| > N_\epsilon$ and $(k_n, m_n) \in C(\alpha, \delta_\epsilon)$, then*

$$\left| \frac{f(k_n, m_n)}{|k_n| + |m_n|} - \frac{c_1}{1 + |\alpha|} h - \frac{c_2 |\alpha|}{1 + |\alpha|} v \right| < \epsilon, \tag{1}$$

*for appropriate $c_1$ and $c_2$. When $\alpha = \pm\infty$, $\delta_\epsilon$ should be replaced by a constant $M_\epsilon$ and the condition $(k_n, m_n) \in C(\alpha, \delta_\epsilon)$ should be replaced by $|m_n/k_n| > M_\epsilon$. Moreover, similar statements are valid for all other cases of Lemma 2.1.*

*Proof.* Draw a uniform $(0, 1)$ distributed random variable $U$ and consider the line $y = \alpha X + U$. Let $y_n$ be the (random) point on the vertical line $\{(x, y) : x = n\}$ closest to this line. It is not hard to see that $(f(y_0), f(y_1), \ldots)$ forms a random walk with (dependent) stationary increments. We write $x_n = (k_n, m_n)$ and write $y_{j(n)}$ for the (or a) vertex among $(y_0, y_1, \ldots)$ which is closest to $x_n$. We then have

$$\frac{f(x_n)}{\|x_n\|} = \left( \frac{f(y_{j(n)})}{\|y_{j(n)}\|} + \frac{f(x_n) - f(y_{j(n)})}{\|y_{j(n)}\|} \right) \frac{\|y_{j(n)}\|}{\|x_n\|}.$$

Since $(f(y_0), f(y_1), \ldots)$ has stationary increments, the ergodic theorem tells us that $f(y_n)/\|y_n\|$ converges a.s., and it then follows from the corresponding convergence in probability in Lemma 2.1 that this a.s. limit must be the same limit as in Lemma 2.1. Therefore, if $k_n$ and $m_n$ are large enough and $|m_n/k_n - \alpha|$ is small enough, then $j(n)$ is large and therefore $f(y_{j(n)})/\|y_{j(n)}\|$ is close to the correct limit in Lemma 2.1. At the same time, the term $\|y_{j(n)}\|/\|x_n\|$ is close to 1 by construction. Finally, the norm of the vector $(f(x_n) - f(y_{j(n)}))/\|y_{j(n)}\|$ is bounded above by

$$\frac{\|x_n - y_{j(n)}\|}{\|y_{j(n)}\|}.$$

This last expression is close to 0 when $k_n$ and $m_n$ are large and $|m_n/k_n - \alpha|$ is small. $\square$

For $\alpha = \alpha_0 := -\frac{h}{v}$, the limit in Lemma 2.1(i) is equal to 0, and we write $C_0(\epsilon)$ for $C(\alpha_0, \epsilon)$.

**Lemma 2.3.** *Let $\epsilon > 0$. With $\mu$-probability one, only finitely many points $z$ outside $C_0(\epsilon)$ have $f(z) = 0$.*

*Proof.* The proof is by contradiction. Suppose infinitely many such $z$ exist. Look at the set $B$ of directions $\beta$ outside $C_0(\epsilon)$ such that for every $\delta > 0$, the cone $C(\beta, \delta)$ contains infinitely many $z$ with $f(z) = 0$. The set $B$ is closed and invariant under translations and therefore $\bar{\beta} := \sup\{\beta : \beta \in B\}$ is well defined and an a.s. constant. According to Lemma 2.2, for every $\gamma > 0$, we now have infinitely many $z$ with $f(z) = 0$ for which

$$\left| \frac{f(z)}{\|z\|} - \frac{c_1 \bar{\beta}}{1 + |\bar{\beta}|} h - \frac{c_2 \bar{\beta} |\bar{\beta}|}{1 + |\bar{\beta}|} v \right| < \gamma.$$

But since $f(z) = 0$ and $\gamma$ is arbitrary, this implies that

$$\frac{c_1 \bar{\beta}}{1 + |\bar{\beta}|} h + \frac{c_2 \bar{\beta} |\bar{\beta}|}{1 + |\bar{\beta}|} v = 0,$$

which implies that $\bar{\beta} = \alpha_0$, a contradiction. □

## 3. Percolation

As mentioned in the introduction, we will concentrate on percolation of edges labelled with 0. The *cluster* $C(z)$ of the vertex $z$ is the set of vertices that can be reached from $z$ by travelling over 0-labelled edges only. We are interested in the question whether or not infinite clusters exist and if so, how many.

A subset $S$ of $\mathbb{Z}^2$ is said to have *density* $r$ if for each sequence $R_1 \subseteq R_2 \subseteq \cdots$ of rectangles in $\mathbb{Z}^2$ with $\cup_n R_n = \mathbb{Z}^2$, it is the case that

$$\lim_{n \to \infty} \frac{\#(S \cap R_n)}{\#(R_n)} = r,$$

where $\#(\cdot)$ denotes cardinality. Burton and Keane (1991) showed that for every stationary percolation process $\mathbb{Z}^2$, a.s. all clusters have density. In addition, they showed that either all infinite clusters have positive density a.s., or all infinite clusters have zero density a.s.

**Lemma 3.1.** *For any cocycle measure $\mu$, all clusters have zero density a.s.*

*Proof.* If infinite clusters exist with positive probability, then the probability that the cluster of the origin is infinite must be positive. All elements $z$ in this cluster have $f(z) = 0$. But according to Lemma 2.3, this implies that for all $\epsilon > 0$, up to a finite number of vertices, the whole cluster is contained in $C_0(\epsilon)$. Since the density of $C_0(\epsilon)$ goes to zero with $\epsilon$ tending to zero, we find that the cluster of the origin has density zero a.s. The result now follows from the result of Burton and Keane just mentioned. □

**Lemma 3.2.** *If infinite clusters exist with positive $\mu$-probability (and hence with probability one according to the ergodicity of $\mu$), then there are infinitely many infinite clusters $\mu$-a.s.*

*Proof.* There are at least two quick proofs of this fact. (1): If infinite clusters exists, then the set of vertices $z$ for which $C(z)$ is infinite has positive density a.s. Since each single cluster has zero density, the only conclusion is that there are infinitely many clusters. (Note that density is not countably additive!) (2): If there are infinite clusters a.s., then there are infinitely many points $z$ on the $x$-axis for which $C(z)$ is infinite. It is clear that $f(z)$ takes infinitely many values among these points $z$. However, when $f(z) \neq f(z')$, then $C(z) \cap C(z') = \emptyset$, since $f(z)$ is obviously constant on a cluster. Hence infinitely many infinite clusters must exist. □

The following result could be stated in higher generality. Strictly speaking, it follows from a general result like Lemma 2.3 in Häggström and Meester (1996), but for this particular situation an independent simple proof is possible. See also Meester (1999) for related results.

**Lemma 3.3.** *For any cocycle measure $\mu$, with probability one every infinite cluster $C$ satisfies*

$$\sup\{y\,;\, (x,y) \in C\} = \infty,$$

*and*

$$\inf\{y\,;\, (x,y) \in C\} = -\infty.$$

*Similar statements are true for the horizontal direction.*

*Proof.* Since $C$ is contained (up to finitely many points) in every $C_0(\epsilon)$, the assumption that $C$ is infinite implies that $C$ is unbounded in at least one of the two vertical directions. We now assume (wlog) that with positive probability (and therefore with probability 1) a cluster $C$ exists for which $\sup\{y\,;\, (x,y) \in C\} = \infty$ and $\inf\{y\,;\, (x,y) \in C\} > -\infty$. According to Lemma 2.3, the intersecion of $C$ with a horizontal line contains at most finitely many points, and therefore there is a leftmost point $C_\ell(n)$ at every level $y = n$, for $n$ large enough. The collection $\{C_\ell(n)\}$, where $C$ ranges over all upwards unbounded clusters of the halfspace $\{y \geq n\}$, has a well defined (one-dimensional horizontal) density $d_n$. It is clear from the construction that $d_{n+1} \geq d_n$, since the restriction to $\{y \geq n+1\}$ of every upwards unbounded cluster of $\{y \geq n\}$ contains at least one upwards unbounded cluster in $\{y \geq n+1\}$. At the same time, $(d_n)$ forms a stationary sequence. We conclude that $d_n$ is constant a.s. On the other hand, note that our assumption implies that any given vertex $z$ is the left-lowest point of an upwards unbounded cluster with positive probability. Therefore the line $y = n+1$ contains a positive density of such points. Clearly, these points are 'new' in the sense that they are not in previous clusters. From this we see that $d_{n+1} > d_n$, the required contradiction. □

Next, we show that percolation occurs in a directed sense:

**Lemma 3.4.** *Every infinite cluster $C$ contains a strictly northwest-southeast directed bi-infinite path. More precisely, the left boundary of $C$ forms such a path.*

*Proof.* Define, for every $n$, the leftmost vertical edge in $C$ between $\{y = n\}$ and $\{y = n+1\}$ by $e_n$. Connect, for all $n$, the upper endpoint of $e_n$ with the lower endpoint of $e_{n+1}$ through the horizontal edges in between them (if necessary). I claim that the union of these vertical and horizontal edges is the required path. To see this, note that the upper endpoint of $e_n$ is connected by a path of 0-edges to the upper endpoint of $e_{n+1}$ since they belong to the same cluster $C$. Since we can also travel from the upper endpoint of $e_n$ to the upper endpoint of $e_{n+1}$ by first travelling horizontally to the lower endpoint of $e_{n+1}$ and the last step via $e_{n+1}$, all these last travelled horizontal edges must have zero labels. Finally, $e_{n+1}$ cannot be strictly at the right of $e_n$, since then the path constructed by traveling vertically from the upper endpoint of $e_n$ and then horizontally to the upper endpoint of $e_{n+1}$ would consist of zero labels only, which contradicts the definition of $e_{n+1}$. □

Finally, we show that 'dead ends' are impossible:

**Lemma 3.5.** *Consider the following event: there is a directed path from the origin going down-right, which is completely labelled 0, and the two edges $((-1,0),(0,0))$ and $((0,1),(0,0))$ are both labelled 1. This event has probability 0.*

*Proof.* According to the previous lemma, the left boundary of the 0-cluster of the origin forms a bi-infinite directed path $\pi$. This bi-infinite path crosses the $x$-axis to the left of the origin. The cluster of the origin must contain a connection between $\pi$ and the directed path going down from the origin. Now label all edges which are forced to be zero, given this connection. It is easy to see that these 0's run into conflict with one of the two designated edges having label 1. □

By now we have a fairly precise and specific description of the geometry of infinite clusters in cocycle measures, if they exist: there are in that case infinitely many such clusters, essentially contained in a cone in the $\alpha_0$-direction, and bounded at the left by a bi-infinite directed path. This description is so specific that it makes it easy in many case to rule out percolation almost immediately. On the other hand, one might suspect that this specific description makes it almost impossible for 'natural' cocycle measures to percolate, and that a percolating cocycle measure must be more or less constructed for that purpose. It seems hard to formulate a general property which excludes percolation. One is tempted to try to connect certain ergodic-theoretical mixing properties with percolation here, since the above description of percolation clusters seems highly non-mixing. However, we will see that mixing cocycle measures that percolate can be constructed; they can even have trivial full tail.

## 4. An example

The following example is discussed to some extent in Burton, Dajani and Meester (1998). Here we shall discuss the example in detail. Choose $0 < p < 1$ and let $q = 1 - p$. Label all edges of the $x$-axis 0 with probability $q$ and 1 with probability $p$, independently of each other. For the $y$-axis we do the same with interchanged probabilities. Now denote the square $[n, n+1] \times [0, 1]$ by $W_n$, and denote the lower and upper edge of $W_n$ by $e_n$ and $f_n$ respectively. The labelling procedure is at follows: first label the remaining edges of $W_1$; if there are two possibilities for doing this we choose one of them with equal probabilities. At this point, the lower and left edge of $W_2$ are labelled, and we next complete the labelling of $W_2$, noting again that if there are two ways to do this, we choose one of them with equal probabilities. This procedure is continued and gives all labels in the strip $[0, \infty) \times [0, 1]$. Then we move one unit upwards, and complete in a similar fashion the labels in the strip $[0, \infty) \times [1, 2]$. (Of course, if you want to carry out this labelling, you never actually finish any strip. Instead, you start at some moment with the second strip which can be labelled as far as the current labelling of the first strip allows, etc.)

This procedure yields a random labelling of all edges in the first quadrant. Using for instance Kolmogorov's consistency condition, we can extend this to a cocycle measure on the labels in the whole plane.

**Lemma 4.1.** *The procedure described above yields a stationary and mixing measure $\mu$. In particular, the labelling is ergodic.*

*Proof.* If we can show that the labelling of the edges $f_n$ has the same distribution as the labelling of the edges $e_n$, then we have shown that the labelling in the quadrant $[0, \infty) \times [1, \infty)$ has the same distribution as the labelling in $[0, \infty) \times [0, \infty)$ and we can use a similar argument for vertical lines plus induction to finish the argument. Therefore we only need to show that the labelling of the edges $f_n$ is i.i.d. with the correct marginals. To do this properly, consider the labels of the edges of $W_n$. There are six possible labellings of the edges of $W_n$. Four of these are such that $e_n$ and $f_n$ have the same label. The exceptional labellings are (starting at the lower left vertex and moving clockwise) 0110 and 1001. Denote the labelling of the edges of $W_n$ by $L_n$. Then it is not hard to see that $L_n$ is a Markov chain on the state space $\{0110, 1001, 1010, 0101, 1111, 0000\}$. Take the transition matrix $P$ of $L_n$, interchange the rows and the columns corresponding to 0110 and 1001 to obtain $P'$, and consider the *backward* Markov chain corresponding to $P'$, denoted by $M_n$. An easy calculation then shows that $L_n$ and $M_n$ have the same transition matrix and that they are both in stationarity. But now note that $M_n$ represents the right-to-left labelling of the strip $[0, \infty) \times [0, 1]$. This means that the distribution of the random vector $(f_{k_1}, f_{k_2}, \ldots, f_{k_n})$ is the same as the distribution of $(e_{k_n}, e_{k_{n-1}}, \ldots, e_{k_1})$. The last vector has independent marginals, hence so has the first, and we are done.

Next we show that $\mu$ is mixing. For this, consider finite-dimensional cylinder events $A$ and $B$, i.e. $A$ and $B$ only depend on edges in the box $B_n = [0, n]^2$. Denote by $M_k$ the labelling of all the edges in the box $B_n + (k, 0)$. It is easy to check that $M_k$ is a mixing Markov chain. This implies that for the events $A$ and $B$, we have

$$\mu(A \cap T_{(k,0)}^{-1}(B)) \to \mu(A)\mu(B),$$

where $T_z$ denotes translation over the vector $z$. This shows that $\mu$ is mixing in the horizontal direction. For the vertical direction, we consider the Markov chain associated with the labellings of the boxes $B_n + (0, k)$, $k = 0, 1, \ldots$, and repeat the argument. □

**Theorem 4.2.** *The measure $\mu$ described in this subsection does not percolate a.s.*

Before we give a proof, we need to look at the construction of $\mu$. The above definition of $\mu$ is simple, but has the disadvantage that we need to appeal to Kolmogorov's consistency theorem to define it on the whole plane. In this sense, the definition is not constructive. There is an alternative way of defining $\mu$ that is constructive, and that will be quite useful. The first step towards this construction is the following lemma.

**Lemma 4.3.** *Let $\pi$ be a bi-infinite path $(\ldots, z_{-1}, z_0, z_1, \ldots)$, where $z_i = (z_{i_1}, z_{i_2})$, and with the property that $z_{k_1}$ is non-decreasing in $k$, and $z_{k_2}$ is non-increasing in $k$. Denote the edge $(z_i, z_{i+1})$ by $e_i$. Then the labels $(\ldots, c(e_{-1}), c(e_0), c(e_1), \ldots)$ form an independent sequence.*

*Proof.* Independence is defined in terms of finite collection of edges, so by stationarity we need only look at finite paths with these monotonicity properties which travel from the $y$-axis to the $x$-axis. That is, we only work in the first quadrant for now.

Denote such a path by $(z_0, \ldots, z_k)$, where $z_0$ is on the $y$-axis, and $z_k$ is on the $x$-axis. We may assume that $z_0$ and $z_k$ are the only points of the path on the coordinate axes. Again denoting the edge $(z_i, z_{i+1})$ by $e_i$, We claim that the last edge $c(e_{k-1})$ is independent of the collection $\{c(e_0), \ldots, c(e_{k-2})\}$. To see this, first assume that $e_{k-2}$ is a horizontal edge. (Note that $e_{k_1}$ is always vertical by assumption.) Then by considering the reversed Markov chain in the proof of Lemma 4.1, it follows immediately that $c(e_{k-1})$ is independent of $c(e_{k-2})$ and is also independent of the labels of all horizontal edges to the left of $e_{k-2}$. But all labels $\{c(e_1), \ldots, c(e_{k-2})\}$ are measurable with respect to these last labels together with some independent labels on the $y$-axis and some independent choices when appropriate. This proves the claim. If $e_{k-2}$ is vertical, we walk back along the path until the first horizontal edge, and repeat the argument with the Markov chain corresponding to the strip with the appropriate width. The lemma now follows with induction in the obvious way. □

The last lemma tells us how to construct a labelling of the whole plane in a constructive manner. We first take a bi-infinite northwest-southest directed path $\pi$ that is unbounded in all directions. Label the edges of this path in an independent fashion, with the correct one-dimensional marginals. Next start 'filling the plane' in a way similar to the original construction. Above the path, we can do what we did before, and label strips from left to right; below the path we use the backwards Markov chain mentioned in the proof of Lemma 4.1, and we label strips from right to left. It is clear that the labelling obtained this way has the correct distribution: just note that all finite-dimensional distributions are correct. Here we need Lemma 4.3 of course, to make sure we start off with the correct distribution on our path $\pi$. Note that the labellings above and below $\pi$ are conditionally independent, given the labelling of $\pi$ itself. Also note that the sigma fields generated by $\mathcal{E}_1(\pi) := \{e : e$ is both below $\pi$ and to the right of the line $x = 0$ (inclusive)$\}$ and $\mathcal{E}_2 := \{e : e$ is both above $\pi$ and to the left of the line $x = 0$ (inclusive)$\}$ are independent.

*Proof of Theorem 4.2.* We assume $\mu$ does percolate, and show that this leads to a contradiction. Choose a northwest-southeast directed path $\pi = (\ldots, z_{-1}, z_0, z_1, \ldots)$ as follows (recall the definition of $h$ and $v$): $z_0$ is the origin, and for some (possibly large) $M$, $z_{-M}, \ldots, z_M$ are all on the $x$-axis. The vertices $z_{M+1}, z_{M+2}, \ldots$ can now be chosen in such a way that they are all above the line through $z_M$ with direction $-h/2v$; the vertices $z_{-M-1}, z_{-M-2}, \ldots$ can now be chosen in such a way that they are all below the line through $z_{-M}$ with direction $-h/2v$. Label the edges of $\pi$ independently (with the correct marginals). The point of this choice is that according to our geometrical picture obtained in the previous section, for some $M$, there must be a positive probability that the origin is contained in a 0-labelled directed infinite path going down-right, all whose edges are strictly below $\pi$. We call this event $E_1$. It is clear that $E_1$ is measurable with respect to the sigma field generated by $\mathcal{E}_1(\pi)$.

On the other hand, there is a positive probability that the four edges between the origin, $(0, 1)$, $(-1, 0)$ and $(-1, -1)$ are all labelled 1. This event, let's call it $E_2$, is measurable with respect to the sigma field generated by $\mathcal{E}_2(\pi)$. We noted above that these two sigma fields are independent and it follows that $E_1$ and $E_2$ are independent. Hence $P(E_1 \cap E_2) > 0$, but on this event, the origin is a dead end in the sense of Lemma 3.5, a contradiction. □

It is interesting to compare the last theorem with a result of Kesten (1982). He showed that if we label all edges independently, with vertical edges being 0 with probability $p$ and horizontal edges with probability $1 - p$, then the system does not percolate. In fact, the system is critical in the sense that if one increases the probability for either the horizontal or vertical edges by a positive amount, the sytem does in fact percolate.

## 5. A percolating cocycle with trivial full tail

The geometrical picture of infinite clusters looks highly non-mixing. After all, if we look at the realisation below the horizontal line $y = n$, we have a lot of information about the infinite clusters, an this should tell a lot about the realisation in the halfplane $y \geq n + m$, for $m$ large. So it seems that a percolating cocycle measure has long distance dependencies and therefore weak mixing properties. But we shall now see that a percolating cocycle measure can be constructed which has trivial full tail, which is much stronger than being mixing. The construction is based on an exclusion process introduced in Yaguchi (1986) and studied by Hoffman (1998).

We first describe Yaguchi's construction. We shall initially work in the half plane $x \geq 0$, but the measure can of course be extended to a measure on the full plane. Consider the $y$-axis. Each vertex is either blue, red or not coloured (probability comes in later). We next colour the line $x = 1$ as follows. Each coloured vertex $z$ (on the $y$-axis) decides independently with a certain (fixed and constant) probability if it wants to move down one unit. It also checks whether or not the vertex below is not coloured. If both the vertex wants to move, and the vertex below is not coloured, then we colour the vertex $z+(1,-1)$ with the same colour as $z$. Otherwise we colour $z + (1,0)$ with the same colour as $z$. This procedure is repeated when we go from $x = 1$ to $x = 2$, etc. Yaguchi (1986) characterised the stationary measures of the associated $\mathbb{Z}^2$ action, in particular he showed such measures exist. Hoffman (1999) showed that these measures have trivial full tail; in particular they are mixing.

How does this relate to cocycle measures? We shall make a few minor modifications. First, we look at all red points on the $y$-axis that are between two given blue points (with no other blue points in between). Take the top vertex among these red vertices and change its colour into green. When a coloured vertex $z$ causes $z+(1,0)$ to be coloured with the same colour as $z$ we also colour the edge between these two vertices with the same colour. When $z$ causes $z + (1,-1)$ to be coloured, we colour the two edges $(z, z + (0,-1))$ and $(z + (0,-1), z + (1,-1))$ with the same colour. Finally we keep the green edges and 'uncolour' all other edges. A configuration now consists of infinitely many disjoint, bi-infinite strictly northwest-southeast directed green paths. We can transform a realisation to a labelling of the edges that satisfy the cocycle identity as follows: All green edges are labelled 0. All edges that have one endpoint in common with a green edge are label 1. All remaining edges are labelled 0. It is easy to prove that the realisation obtained this way satisfies our cocycle identity. It is obtained as a ergodic-theoretical factor of a mixing process, and therefore also mixing. On the other hand, it percolates along the edges that are coloured green.